\documentclass[a4paper, oneside,11pt]{article}

\usepackage[latin1]{inputenc}
\usepackage[T1]{fontenc}
\usepackage[english]{babel}
\usepackage{amssymb}
\usepackage{amsmath}
\usepackage{amsthm}
\usepackage{graphicx}

\usepackage[all]{xy}
\usepackage{verbatim}
\usepackage{amsmath}
\usepackage{amsthm}
\usepackage{graphicx}
\usepackage{amssymb}
\usepackage{epstopdf}
\usepackage{url}
\usepackage{amsfonts}
\usepackage{todonotes}
\newcommand{\mc}{\mathcal}

\newcommand{\br}{\mathbb{R}}

\newcommand{\bt}{\mathbb{T}}

\renewcommand{\)}{\right)}

\newtheorem{thm}{Theorem}
\newtheorem{lem}[thm]{Lemma}

\newtheorem{prop}[thm]{Proposition}

\newtheorem{remark}[thm]{Remark}

\def\be{\begin{equation}}
\def\ee{\end{equation}}
\def\bea{\begin{eqnarray}}
\def\eea{\end{eqnarray}}

\usepackage{enumerate}

\usepackage{bbm}

\numberwithin{thm}{section}
\numberwithin{equation}{section}

\renewcommand{\div}{\operatorname{div}}

\usepackage{hyperref}
\usepackage{mathtools}
\mathtoolsset{showonlyrefs}

\newcommand*\di{\mathop{}\!\mathrm{d}}

\newcommand{\RR}{\mathbb{R}}

\newcommand{\TT}{\mathbb{T}}

\newcommand{\ds}{\displaystyle}

\title{Recent developments on the well-posedness theory \\ for Vlasov-type equations}

\author{
Megan Griffin-Pickering 
  \thanks{Durham University, Department of Mathematical Sciences, Lower Mountjoy, Stockton Road, Durham DH1 3LE, UK.  Email: \textsf{megan.k.griffin-pickering@durham.ac.uk}}
  \and
Mikaela Iacobelli
  \thanks{ETH Z\"urich, R\"amistrasse 101, 8092 Z\"urich, Switzerland. Email: \textsf{mikaela.iacobelli@math.ethz.ch}}
}

\begin{document}

\maketitle

\abstract{In these notes we summarise some recent developments on the existence and uniqueness theory for Vlasov-type equations, both on the torus and on the whole space. }

\section{An introduction to Vlasov-type equations in plasma physics}

In this note, we discuss some recent results concerning a class of PDEs used in the modelling of plasma.
Plasma is a state of matter abundant in the universe. It can be found in stars, the solar wind and the interstellar medium, and is therefore widely studied in astrophysics, as well as in many other contexts. For example, a major terrestrial application is in nuclear fusion research. 
For this reason, mathematical modelling of plasma is of interest,
with different types of plasma models being suitable for different contexts.

A plasma consists of an ionised gas. It forms when an electrically neutral gas is subjected to high temperatures or a strong electromagnetic field, which causes the gas particles to dissociate into charged particles. These charged particles then interact through electromagnetic forces. The relatively long range nature of these interactions results in a collective behaviour distinct from that expected from a neutral gas.

In this article, we will discuss the well-posedness of a certain class of PDE models for plasma. We will consider equations of Vlasov type, which describe particle systems with mean field interactions.

\subsection{The Vlasov-Poisson system: the electrons' view-point}

The ionisation process in the formation of a plasma produces two types of charged particle: positively charged ions and negatively charged electrons. It also generally contains neutral species, since not all of the particles of the original neutral gas will dissociate. However, typically the interactions with the neutral species are weak in comparison to the interactions of the charged species. For the purposes of these notes, we will neglect interactions with the neutral particles and concentrate on the modelling of the charged particles. 

In fact, it is usual to make an assumption which decouples the dynamics of the two species. It is possible to do this because the mass of an electron is much smaller than the mass of an ion. 
The result is a separation between the timescales on which each species evolves: in short, the ions typically move much more slowly than the electrons.
When modelling the electrons, it is thus common to assume that the ions are stationary over the time interval of observation.

The Vlasov-Poisson system is a well-known kinetic equation describing this situation. This equation was proposed by Jeans \cite{Jeans} as a model for galaxies. Its use in the plasma context dates back to the work of Vlasov \cite{Vlasov1}. 
The most commonly known version of the system models the electrons in the plasma.
The electrons are described by a density function $f = f(t,x,v)$, which is the unknown in the following system of equations: 
\be \label{eq:VP-physical-bg}
(VP) : =  \begin{cases}
\partial_t f + v \cdot \nabla_x f + \frac{q_e}{m_e} E \cdot \nabla_v f = 0, \\
\nabla_x \times E = 0, \\
\epsilon_0 \nabla_x \cdot E = R_i + q_e \rho_f , \\ \ds
\rho_f(t,x) : = \int_{\RR^d} f(t,x,v) \di v, \\
f \vert_{t=0} = f_0 \geq 0 .
\end{cases}
\ee
Here $q_e$ is the charge on each electron, $m_e$ is the mass of an electron and $\epsilon_0$ is the electric permittivity. $R_i : \RR^d \to \RR_+$ is the charge density contributed by the ions, which is independent of time since they are assumed to be stationary. The electrons experience a force $q_e E$, where $E$ is the electric field induced by the plasma itself. This is found from the Gauss law
\be
\nabla_x \times E = 0, \quad \epsilon_0 \, \nabla_x \cdot E = R_i + q_e \rho_f,
\ee
which arises as an electrostatic approximation of the full Maxwell equations.

The system \eqref{eq:VP-physical-bg} expresses the fact that each electron in the plasma feels the influence of the other particles in the plasma in an averaged sense, through the electric field $E$ induced collectively by the whole plasma. This is a long-range interaction between particles. In particular, this equation does not account for collisions between particles of any species.

The Vlasov-Poisson system as written in equation \eqref{eq:VP-physical-bg} does not yet include a boundary condition.
In this note we focus on two cases: either the periodic case
where the spatial variable $x$ lies in the $d$-dimensional flat torus $\TT^d$
and the velocity variable $v$ lies in the whole Euclidean space $\RR^d$, or the whole space case where both $x$ and $v$ range over $\RR^d$.
We will use the notation $\mc{X}$ to denote the spatial domain, either $\TT^d$ or $\RR^d$ as appropriate, so that throughout this note we have $(x,v) \in \mc{X} \times \RR^d$.

It is common to restrict in particular to the case where the background ion density $R_i$ is spatially uniform. 
In the case of the torus, $\mc{X} = \TT^d$, this results in the system
\be \label{eq:VP-physical-torus}
(VP) : =  \begin{cases}
\partial_t f + v \cdot \nabla_x f + \frac{q_e}{m_e} E \cdot \nabla_v f = 0, \\
\nabla_x \times E = 0, \\
\epsilon_0 \nabla_x \cdot E = q_e \left (\rho_f - \ds \int_{\TT^d \times \RR^d} f \di x \di v \right ) , \\ \ds
f \vert_{t=0} = f_0 \geq 0 .
\end{cases}
\ee
The ion charge density is chosen to be 
\be
R_i \equiv - q_e \int_{\TT^d \times \RR^d} f \di x \di v
\ee
so that the system is globally neutral. 
This is required from the point of view of the physics under consideration due to the conservation of charge, since the plasma forms from an electrically neutral gas.
Note that any solution $f$ of \eqref{eq:VP-physical-torus} satisfies a transport equation with a divergence free vector field, which implies that the mass of $f$ is conserved by the evolution. Thus in fact
\be
R_i \equiv - q_e \int_{\TT^d \times \RR^d} f_0 \di x \di v .
\ee
In mathematical treatments, it is common to see \eqref{eq:VP-physical-torus} written in the rescaled form
\be \label{eq:VP}
(VP) : =  \begin{cases}
\partial_t f + v \cdot \nabla_x f + E \cdot \nabla_v f = 0, \\
\nabla_x \times E = 0, \\
 \nabla_ x \cdot E = \rho_f - 1, \\ \ds
f \vert_{t=0} = f_0 \geq 0 , \, \int_{\TT^d \times \RR^d} f_0 \di x \di v = 1.
\end{cases}
\ee

In the whole space case $\mc{X} = \RR^d$, one often considers a vanishing background, in order to have a system with finite mass. This results in the system
\be \label{eq:VP-whole}
\begin{cases}
\partial_t f + v \cdot \nabla_x f + E \cdot \nabla_v f = 0, \\
\nabla_x \times E = 0, \\
 \nabla_ x \cdot E = \rho_f , \\ \ds
f \vert_{t=0} = f_0 \geq 0 , \, \int_{\RR^d \times \RR^d} f_0 \di x \di v = 1.
\end{cases}
\ee

\subsection{The Vlasov-Poisson system with massless electrons: the ions' view-point}

The previous section presented the Vlasov-Poisson system as a model for the electrons in a dilute, unmagnetised, collisionless plasma.
A variant of the Vlasov-Poisson system may be used to model the ions in the plasma instead. 

To derive an appropriate model, once again we make use of the large disparity between the masses of the two species. The resulting separation of timescales allows an approximation in which the two species are modelled separately.
From the point of view of the ions, the electrons have a very small mass and so are very fast moving.
Since the electrons are not stationary, a model of the form \eqref{eq:VP} is not appropriate.

Instead observe that, since the electrons move quickly relative to the ions, the frequency of electron-electron collisions is high in comparison to ion-ion or ion-electron collisions. Electron-electron collisions are expected to be relevant on the typical timescale of evolution of the ions, even while the frequencies of other kinds of collisions remain negligible.
The expected effect of the electron-electron collisions is to drive the electron distribution towards its equilibrium configuration.
In ion models it is therefore common in physics literature to assume that the electrons are close to thermal equilibrium. 

In the limit of \textbf{massless electrons}, the ratio between the masses of the electrons and ions, $m_e/m_i$, tends to zero. Here $m_e$ is the mass of an electron and $m_i$ is the mass of an ion. 
In the limiting regime, it is assumed that the electrons instantaneously assume the equilibrium distribution.
This approximation is often made in the physics literature, motivated by the fact that $m_e/m_i$ is close to zero in applications.

\subsubsection{The Maxwell-Boltzmann Law for Electrons}

The equilibrium distribution can be identified by studying the equation for the evolution of electrons. 
Let the ion density $\rho[f_i]$ be fixed, and assume that all ions carry the same charge $q_i$. We have discussed that a possible model for the evolution of the electron density is the Vlasov-Poisson system \eqref{eq:VP-physical-bg}.
However, the Vlasov-Poisson system is a collisionless model. 
As discussed above, in the long time regime we consider we expect the effect of electron-electron collisions to be significant.  

Collisions in a plasma are described by the Landau-Coulomb operator $Q_L$ \cite[Chapter 4]{Lifshitz-Pitaevskii}, which is an integral operator defined as follows: for a given function $g = g(v) : \RR^d \to \RR$,
\be
Q_L(g) : = \nabla_v \cdot \int_{\RR^d} a(v - v_*) : \left [ g(v_*) \nabla_v g(v) - g(v) \nabla_v g(v_*) \right ] \di v_* .
\ee
The tensor $a$ is defined by
\be
a(z) = \frac{|z|^2 - z \otimes z}{|z|^3} .
\ee
We add this term to the Vlasov-Poisson system to model a plasma with collisions. This results in the following model for the electron density $f_e$:
\be \label{eq:VP-electrons}
\begin{cases} \ds
\partial_t f_e + v \cdot \nabla_x f_e + \frac{q_e}{m_e} E \cdot \nabla_v f_e = \frac{C_e}{m_e^2} Q_L(f_e), \\
\nabla_x \times E = 0, \quad
\epsilon_0 \nabla_x \cdot E = q_i \rho[f_i] + q_e \rho[f_e] .
\end{cases}
\ee
Here $C_e$ is a constant depending on physical quantities such as the electron charge $q_e$ and number density $n_e$, but \textit{not} on the electron mass $m_e$. For the derivation of the scaling $C_e/m_e^2$ in front of the Landau-Coulomb operator, see Bellan \cite[Chapter 13, Equation (13.46)]{Bellan}.

Consider the rescaling
\be
F_e(t,x,v) = m_e^{-\frac{d}{2}} f_e \left ( t, x, \frac{v}{\sqrt{m_e}}\right) .
\ee
Notice that this scaling preserves the macroscopic density: $\rho[F_e] = \rho[f_e]$.
Then $F_e$ satisfies
\be \label{eq:VP-electrons-mass}
 \left \{
\begin{array}{c} \ds
\sqrt{m_e} \partial_t F_e + v  \cdot \nabla_x F_e + q_e E \cdot \nabla_v F_e = C_e Q_L(F_e), \\ \ds
\nabla \times E = 0, \quad
\epsilon_0 \nabla \cdot E = q_i \rho[f_i] + q_e \rho[F_e] .
\end{array} \right.
\ee
We assume that $F_e$ converges to a stationary distribution $\bar f_e = \bar f_e(x,v)$ as $m_e$ tends to zero, and focus on formally identifying $\bar f_e$. 

To identify the possible forms of $ \bar f_e$, we consider the entropy functional
\be 
H[f] : = \int_{\mc{X} \times \RR^d} f \log{f} \di x \di v . 
\ee 
For a solution $F_e$ of Equation~\eqref{eq:VP-electrons-mass},
\be
\frac{\di}{\di t} H[F_e] = m_e^{-1/2} \int_{\mc{X} \times \RR^d} (1 + \log{F_e}) \left [ - \div_{x,v} \left ( (v, E) F_e \right) + Q_L(F_e) \right ] \di x \di v.
\ee
Integrating by parts formally, the transport term vanishes:
\begin{align}
- \int_{\mc{X} \times \RR^d} (1 + \log{F_e}) \div_{x,v} \left ( (v, E) F_e \right) \di x \di v& = \int_{\mc{X} \times \RR^d}  (v, E) \cdot \nabla_{x,v} F_e \di x \di v \\
& = \int_{\mc{X} \times \RR^d} \div_{x,v} \left ( (v, E) F_e \right) \di x \di v = 0 .
\end{align}
Thus
\be
\frac{\di}{\di t} H[F_e] = m_e^{-1/2} \int_{\mc{X} \times \RR^d} (1 + \log{F_e})  Q_L(F_e) \di x \di v.
\ee
By substituting the definition of $Q_L$, one can calculate formally (see \cite{Desvillettes-Villani2000}) that
\begin{multline} \label{def:dissipation}
\frac{\di}{\di t} H[F_e] =\\ - \frac{C}{\sqrt{m_e}} \int_{\mc{X} \times \RR^d} \frac{ 1 }{|v-v_*|} \left \lvert  P_{(v-v_*)^{\perp}}  \left [ \nabla_v \sqrt{F_e(x,v) F_e(x,v_*)} - \nabla_{v_*} \sqrt{F_e(x,v) F_e(x,v_*)} \right ] \right \rvert^2 \di v_* \di v \di x ,
\end{multline}
where $P_{(v-v_*)^{\perp}}$ denotes the operator giving the orthogonal projection onto the hyperplane perpendicular to $v- v_*$. 
For a stationary solution $\bar f_e$, we must have $\frac{\di}{\di t} H[\bar f_e] = 0$, that is, the functional on the right hand side of \eqref{def:dissipation} must vanish. If $\bar f_e \in L^1$, it follows (see for example \cite[Lemma 3]{Villani1996}) that $\bar f_e$ is a local Maxwellian of the form
\be \label{def:bar-fe}
\bar f_e(x,v) = \rho_e(x) \left (\pi \beta_e(x)\right )^{d/2} \exp \left [ - \beta_e(x) |v - u_e(x)|^2 \right ] .
\ee 

The electron density $\rho_e$, mean velocity $u_e$ and inverse temperature $\beta_e$ can then be studied using an argument similar to the one given in the proof of \cite[Theorem 1.1]{BGNS18}.
Substituting the form \eqref{def:bar-fe} into equation \eqref{eq:VP-electrons}, we obtain the following identity for all $x$ such that $\rho_e(x) \neq 0$ and all $v \in \RR^d$:
\begin{multline} 
- \nabla_x \beta_e \cdot (v - u_e) |v - u_e|^2 - u_e \cdot \nabla_x \beta_e |v - u_e|^2 + \beta_e (v - u_e)^\top \nabla_x u_e (v - u_e)\\
 + (v - u_e) \cdot \left [ \nabla_x \log{(\rho_e \beta_e^{d/2})} - q_e \beta_e E + u_e \cdot \nabla_x u_e\right] + u_e \cdot \nabla_x \log{(\rho_e \beta_e^{d/2})} = 0 .
\end{multline}
For each fixed $x$, the left hand side is a polynomial in $v - u_e(x)$, whose coefficients must all be equal to zero. 
For example, by looking at the cubic term we see that $\nabla_x \beta_e = 0$ and thus $\beta_e$ must be a constant independent of $x$.  

The quadratic term then gives
\be
v^\top \nabla_x u_e v = 0 \qquad \text{for all } v \in \RR^d,
\ee
which implies that $\nabla_x u_e$ is skew-symmetric. 
On a spatial domain for which a Korn inequality holds, this restricts the class of $u_e$ that can occur. For example, in the case of the torus $\mc{X} = \TT^d$, the fact that the symmetric part of $\nabla_x u_e$ vanishes implies that $u_e$ is constant \cite[Proposition 13]{Desvillettes-Villani2005}.

Finally, from the linear term we obtain that 
\be
 \nabla_x \log{(\rho_e \beta_e^{d/2})} - q_e \beta_e E = 0.
\ee
Since $\nabla_x \times E = 0$, $E$ is a gradient - that is, it can be written as $E = - \nabla U$ for some function $U$.
Then 
\be
 \nabla_x \log{(\rho_e \beta_e^{d/2})} = - q_e \beta_e \nabla_x U .
\ee
From this we deduce that $\rho_e$ should be of the form 
\be
\rho_e(x) = A \exp \left ( - q_e \beta_e U \right ),
\ee
for some constant $A>0$.
This is known as a \textbf{Maxwell-Boltzmann} law. 

In the whole space case $\mc{X} = \RR^d$, we include an additional spatial confinement of the electrons, by adding an additional potential $\Psi$ to the electron dynamics. 
The equivalent of equation \eqref{eq:VP-electrons-mass} is then
\be \label{}
 \left \{
\begin{array}{c} \ds
\sqrt{m_e} \partial_t F_e + v  \cdot \nabla_x F_e + ( q_e E - \nabla \Psi ) \cdot \nabla_v F_e = C_e Q_L(F_e), \\ \ds
\nabla \times E = 0, \quad
\epsilon_0 \nabla \cdot E = q_i \rho[f_i] + q_e \rho[F_e] .
\end{array} \right.
\ee
Repeating the previous argument, we can derive the following limiting distribution in the regime $m_e \to 0$:
\be
\rho_e = A e^{ - \beta_e (q_e U + \Psi) }  = A g e^{- \beta_e q_e U},
\ee
where we let $g = e^{- \beta_e \Psi}$. We assume that the confining potential $\Psi$ grows sufficiently quickly at infinity so that $g \in L^1 \cap L^\infty(\RR^d \times \RR^d)$.

Bardos, Golse, Nguyen and Sentis \cite{BGNS18} studied the problem of rigorously identifying the Maxwell-Boltzmann law as the distribution of electrons in the massless limit. They consider coupled systems of the form
\be \label{eq:VP-coupled}
\begin{cases}
\partial_t f_i + v \cdot \nabla_x f_i + \frac{q_i}{m_i} E \cdot \nabla_v f_i = 0, \\
\partial_t f_e + v \cdot \nabla_x f_e + \frac{q_e}{m_e} E \cdot \nabla_v f_e = C(m_e) Q(f_e), \\
\nabla_x \times E = 0, \quad
\epsilon_0 \nabla_x \cdot E = q_i \rho[f_i] + q_e \rho[f_e] .
\end{cases}
\ee
In the above, $Q$ denotes a collision operator such as a BGK or Boltzmann operator. Under suitable hypotheses on the spatial domain and the collision rate $C(m_e)$, and assuming the existence of sufficiently regular solutions of the coupled system \eqref{eq:VP-coupled}, they derive that, in the limit as $m_e/m_i$ tends to zero, the electrons indeed take on a Maxwell-Boltzmann distribution. Moreover, solutions of the system \eqref{eq:VP-coupled} converge to a solution of a system of a similar form to \eqref{eq:vpme}, but where the electron temperature depends on time and is chosen to respect the conservation of energy. Other works on this topic include, for example, the work of Bouchut and Dolbeault \cite{Bouchut-Dolbeault95} on the long time limit for the Vlasov-Poisson-Fokker-Planck system for one species -- the massless electrons limit can be related to a long time limit since \eqref{eq:VP-electrons-mass} can also be seen as a time rescaling. Herda \cite{Herda16} also considered the massless electron limit in the case with an external magnetic field. In this case the limiting system is a fluid model for the electrons, coupled with a kinetic model for the ions.

\subsubsection{The Vlasov-Poisson System in the Limit of Massless Electrons}

From equation \eqref{eq:VP-electrons}, we see that the electrostatic potential $U$ induced by a distribution $\rho[f_i]$ of ions with a background of thermalised electrons should satisfy the following semilinear elliptic PDE:
\be \label{eq:semilinear}
- \epsilon_0 \Delta U = q_i \rho[f_i] + A q_e \, g \exp \left ( - \frac{q_e U}{k_B T_e }\right ),
\ee
where in the torus case $\mc{X} = \TT^d$ we let $g \equiv 1$.
The normalising constant $A$ should be chosen so that the system is globally neutral, that is, the total charge is zero:
\be \label{eq:Poisson-SL}
\int_{\mc{X}} q_i \rho[f_i] + A q_e  \exp \left ( - \frac{q_e U}{k_B T_e }\right ) \di x = 0 .
\ee

Indeed, on the torus $\mc{X} = \TT^d$, the Poisson equation
\be
\Delta U = h
\ee
can only be solved if $h$ has total integral zero.
Thus if \eqref{eq:semilinear} has a solution, global neutrality must hold automatically. 
Adjusting the choice of $A$ corresponds to adding a constant to $U$. Thus without loss of generality we choose $A=1$.

Then, the nonlinear equation \eqref{eq:semilinear} replaces the standard Poisson equation for the electrostatic potential in the Vlasov-Poisson system \eqref{eq:VP-physical-bg}. After a suitable normalisation of physical constants, this leads to the following system for the ions:
\be \label{eq:vpme}
(VPME) : = \begin{cases}
\partial_t f + v \cdot \nabla_x f + E \cdot \nabla_v f = 0, \\
E = - \nabla_x U, \\
 \Delta U = e^U - \rho_f , \\ \ds
f \vert_{t=0} = f_0 \geq 0 , \; \int_{\TT^d \times \RR^d} f_0 \di x \di v = 1.
\end{cases}
\ee
This is known as the Vlasov-Poisson system \textbf{with massless electrons}, or \textbf{VPME} system.

In the whole space case $\mc{X} = \RR^d$, we consider two versions of the VPME system, depending on the choice of the constant $A$.
In one case, we let $A=1$. With a suitable choice of dimensionless variables, this results in the following system:
\be \label{eq:vpme-whole-variable}
\begin{cases}
\partial_t f + v \cdot \nabla_x f + E \cdot \nabla_v f = 0, \\
E = - \nabla_x U, \\
 \Delta U = g e^U - \rho_f , \\ \ds
f \vert_{t=0} = f_0 \geq 0 , \; \int_{\RR^d \times \RR^d} f_0 \di x \di v = 1.
\end{cases}
\ee
This system is structurally similar to the torus case \eqref{eq:vpme} considered above. Note however that in this model the system is not necessarily globally neutral.
In order to enforce global neutrality, we can instead choose $A$ to be a normalising constant
\be
A = \frac{1}{\int_{\RR^d} g e^U \di x }.
\ee
Thus we obtain the following system:
\be \label{eq:vpme-whole-fixed}
\begin{cases}
\partial_t f + v \cdot \nabla_x f + E \cdot \nabla_v f = 0, \\
E = - \nabla_x U, \\
 \Delta U = \frac{g e^U}{\int_{\RR^d} g e^U \di x } - \rho_f , \\ \ds
f \vert_{t=0} = f_0 \geq 0 , \; \int_{\RR^d \times \RR^d} f_0 \di x \di v = 1.
\end{cases}
\ee

The VPME system has been used in the physics literature in, for instance, numerical studies of the formation of ion-acoustic shocks \cite{Mason71, SCM} and the development of phase-space vortices behind such shocks \cite{BPLT1991}, as well as in studies of the expansion of plasma into vacuum \cite{Medvedev2011}. A physically oriented introduction to the model \eqref{eq:vpme} may be found in \cite{Gurevich-Pitaevsky75}. 

In \cite{IGP-WP}, we consider the problem of proving well-posedness for the VPME system \eqref{eq:vpme} under reasonable conditions on the initial datum $f_0$. The well-posedness of the systems \eqref{eq:vpme-whole-variable} and \eqref{eq:vpme-whole-fixed} is considered in a forthcoming paper.

\section{Well-posedness for Vlasov equations with smooth interactions}

The Vlasov-Poisson system is an example of a more general class of nonlinear scalar transport equations known as \textbf{Vlasov equations}. A Vlasov equation takes the following form:
\be \label{eq:Vlasov}
\begin{cases}
\partial_t f + v \cdot \nabla_x f + F[f] \cdot \nabla_v f = 0 , \\
F[f](t,x) = - \nabla_x W \ast \rho_f, \\ \ds
\rho_f(t,x) = \int_{\RR^d} f(t,x,v) \di v , \\
f(0,x,v) = f_0(x,v) \geq 0.
\end{cases}
\ee
The system \eqref{eq:Vlasov} is a mean field model for a system of interacting particles with binary interactions described by a pair potential $W: \mc{X} \to \RR$. 
The electron Vlasov-Poisson systems \eqref{eq:VP}, \eqref{eq:VP-whole} can be seen to be of the form \eqref{eq:Vlasov} by choosing $W$ to be the Green's function of the Laplacian on $\mc{X}$. By this we mean that $G$ is a function satisfying the relation
\be \label{eq:fundamental-solution}
- \Delta G = \delta_0 - 1 \quad \text{for } \mc{X} =  \TT^d , \quad \text{ or } \quad - \Delta G = \delta_0 \quad  \text{for }  \mc{X} =  \RR^d .
\ee
The function $U = G \ast ( \rho_f - 1 )$ is a solution of the Poisson equation, respectively
\be
- \Delta U = \rho_f - 1 \quad \text{on} \; \TT^d \quad \text{ or } \quad - \Delta U = \rho_f \quad  \text{on} \;  \RR^d .
\ee
Thus the Vlasov-Poisson systems \eqref{eq:VP}, \eqref{eq:VP-whole} are of the form \eqref{eq:Vlasov}.

The available well-posedness theory for the system \eqref{eq:Vlasov} depends on the choice of the interaction potential $W$, and in particular on the regularity of the force $- \nabla W$. For example, if $\nabla W$ is a Lipschitz function, then the system \eqref{eq:Vlasov} is well-posed in the class $C \left ( [0,\infty); \mc{M}_+(\mc{X} \times \RR^d) \right )$ - the space of continuous paths taking values in the space $\mc{M}_+(\mc{X} \times \RR^d)$ of finite measures on $\mc{X} \times \RR^d$ equipped with the topology of weak convergence of measures. 
This case was considered for example by Braun and Hepp \cite{Braun-Hepp} and by Dobrushin \cite{Dobrushin}.

A path $f \in C \left ( [0,\infty); \mc{M}_+ (\mc{X} \times \RR^d) \right )$ is a weak solution of the Vlasov equation \eqref{eq:Vlasov} if, for all test functions $\phi \in C^1_c \left ( [0,\infty) \times \mc{X} \times \RR^d \right )$,
\be \label{eq:Vlasov-weak}
\int_0^\infty \int_{\mc{X} \times \RR^{d}} \Big[ \partial_t \phi + v \cdot \nabla_x \phi - \left ( \nabla_x  W \ast_x \rho_f \right ) \cdot \nabla_v \phi \Big ] f(t, \di x, \di v) \di t + \int_{\mc{X} \times \RR^{d}} \phi(0, x,v) f_0(\di x, \di v) = 0 .
\ee
Under the assumption that $\nabla W$ is a Lipschitz function, it is known that weak solutions of the Vlasov equation \eqref{eq:Vlasov} exist \cite{Braun-Hepp, Dobrushin} and are unique \cite{Dobrushin}.

\begin{thm} \label{thm:Vlasov-existence}
Assume that $\nabla W : \mc{X} \to \RR^d$ is a Lipschitz function.
Let $f_0$ be a finite non-negative measure with finite first moment:
\be
\int_{\TT^d \times \RR^d} (1 + |x| + |v|) f_0 (\di x \di v) < + \infty .
\ee
Then there exists a unique weak solution $f \in C \left ( [0, + \infty) ; \mc{M}_+(\mc{X} \times \RR^d) \right )$ of the Vlasov equation \eqref{eq:Vlasov}.
\end{thm}

\section{Well-posedness for the Vlasov-Poisson System}

In the case of the Vlasov-Poisson system for electrons \eqref{eq:VP}, the interaction potential $W$ is chosen to be the function $G$ defined by the relation \eqref{eq:fundamental-solution}. The resulting force $K = - \nabla G$ is known as the Coulomb kernel. 
However, $K$ is not a Lipschitz function and so the Vlasov-Poisson system does not satisfy the assumptions of Theorem~\ref{thm:Vlasov-existence}. For example, in the whole space case, $\mc{X} = \RR^d$, $G$ takes the form
\be \label{def:G}
G_{\RR^d}(x) = \begin{cases}
- \frac{1}{2 \pi} \log{|x|}, & d=2, \\
\frac{1}{4 \pi |x|}, & d=3 ,
\end{cases}
\ee
The Coulomb kernel $K_{\RR^d} = - \nabla G_{\RR^d}$ takes the form
\be \label{def:K}
K(x) = \begin{cases}
\frac{x}{2 \pi |x|^2}, & d=2, \\
\frac{x}{4 \pi |x|^3}, & d=3 ,
\end{cases}
\ee
and thus has a singluarity at $x=0$.

On the torus $\mc{X} = \TT^d$, it can be shown that $G_{\TT^d}$ is smooth away from the origin: $G_{\TT^d} \in C^\infty(\TT^d \setminus \{0\})$. Near the singularity it is of the form 
\be
G_{\TT^d} = G_{\RR^d} + G_1,
\ee
where $G_1$ is a $C^\infty$ function. Thus $K_{\TT^d}$ possesses a singularity similar to that of $K_{\RR^d}$.

Consequently, Theorem~\ref{thm:Vlasov-existence} does not apply to the Vlasov-Poisson system. It is not known whether the Vlasov-Poisson system is well-posed in the class of measure solutions. 
However, global well-posedness has been shown for solution classes with greater regularity.

Arsen'ev \cite{Arsenev} introduced a notion of weak solution for the Vlasov-Poisson system \eqref{eq:VP} in dimension $d=3$ and proved the existence of such solutions, globally in time, for initial data $f_0$ belonging to the space $L^1 \cap L^\infty(\RR^6)$.
The boundedness condition $f_0 \in L^\infty(\RR^{6})$ was later relaxed to $f_0 \in L^p(\RR^{6})$, for $p$ sufficiently large, by Horst and Hunze \cite{Horst-Hunze}.

In the case of classical $C^1$ solutions, in the two-dimensional case $d=2$
Ukai and Okabe \cite{Ukai-Okabe} proved global existence for initial data $f_0 \in C^1(\RR^4)$ decaying sufficiently fast at infinity.
In dimension $d=3$, global-in-time solutions were constructed by Pfaffelmoser \cite{Pfaffelmoser} for initial data $f_0 \in C^1_c(\RR^6)$. 
Schaeffer gave a streamlined proof of the same result in \cite{Schaeffer}. 
Horst \cite{Horst1993} extended these results to include non-compactly supported initial data with sufficiently fast decay at infinity. 
The methods of proof for these results are based on an analysis of the characteristic trajectories associated to system \eqref{eq:VP}.
This approach was adapted to the torus by Batt and Rein \cite{Batt-Rein}, who proved the existence of global-in-time classical solutions for \eqref{eq:VP} posed on $\TT^3 \times \RR^3$, for initial data $f_0 \in C^1(\TT^3 \times \RR^3)$ with sufficiently fast decay at infinity.

An alternative approach to the construction of global-in-time solutions in dimension $d=3$ was provided by Lions and Perthame \cite{Lions-Perthame}. Their method is based on proving the propagation of moments. They showed global existence of solutions, provided that the initial datum $f_0 \in L^1 \cap L^\infty(\RR^d \times \RR^d)$ has moments in velocity of sufficiently high order.
However, their strategy is for the whole space case $x \in \RR^d$, and differs from the strategies currently available for the torus.

Pallard \cite{Pallard} then extended the range of moments that could be propagated in the whole space case and showed propagation of moments on the torus $\TT^3$, using a method based on an analysis of trajectories (more similar to \cite{Batt-Rein, Pfaffelmoser, Schaeffer}). Chen and Chen \cite{Chen-Chen} adapted these techniques to further extend the range of moments that could be propagated for the torus case.

Lions and Perthame \cite{Lions-Perthame} proved a uniqueness criterion for their solutions under the additional technical condition that, for all $R, T >0$,
\be
\sup \left\{ |\nabla f_0(y + vt, w)| : |y-x| \leq R, |w-v| \leq R \} \in L^\infty \left ( (0,T) \times \RR^3_x ; L^1 \cap L^2 (\RR^3_v) \right )\right\} .
\ee
Robert \cite{Robert} then proved uniqueness for solutions that are compactly supported in phase space for all time. Subsequently, Loeper \cite{Loeper} proved a uniqueness result which requires only boundedness of the mass density $\rho_f$, and therefore includes the compactly supported case. Loeper's result is based on proving a stability estimate on solutions of the VPME system \eqref{eq:vpme} with bounded density, with respect to their initial data $f_0$ -- in particular, a quantitative estimate in terms of the second order Wasserstein distance $W_2$.
In a similar vein, in the one dimensional case Hauray \cite{Hauray14} proved a weak-strong uniqueness principle, showing that if a bounded density solution exists, then this solution is unique among measure-valued solutions. This result is also based on a Wasserstein stability result.

\section{Well-posedness theory for the Vlasov-Poisson system with massless electrons}

The VPME system for ions is in general less well understood than the Vlasov-Poisson system for electrons, due to the additional nonlinearity in the elliptic equation for the electrostatic potential.
In the case of the well-posedness theory, weak solutions for the VPME system were constructed in dimension $d=3$ in the whole space by Bouchut \cite{Bouchut}, globally in time.
In one dimension, global-in-time weak solutions were constructed by Han-Kwan and Iacobelli \cite{IHK1} for measure data with a first moment. A weak-strong uniqueness principle was also proved for solutions satisfying $\rho_f \in L^\infty_{\text{loc}} \left ([0, + \infty) ; L^\infty(\TT) \right )$: namely, if a solution with this regularity exists, then it is unique among measure solutions. However, a well-posedness theory for strong solutions in higher dimensions remained open.

In the article \cite{IGP-WP}, global well-posedness is proved for the VPME system on the torus in dimension $d=2$ and $d=3$. The main result is stated in the following theorem.

\begin{thm}[Global well-posedness: $\TT^d$] \label{thm:main}
Let $d = 2, 3$. Let the initial datum $f_0 \in L^1 \cap L^\infty(\bt^d \times \br^d)$ be a probability density satisfying
\be
f_0(x,v) \leq \frac{C_0}{1 + |v|^{k_0}}\;\;\mbox{for some}\,\, k_0 > d, \quad \int_{\TT^d \times \RR^d} |v|^{m_0} f_0(x,v) \di x \di v < + \infty\;\;\mbox{for some}\,\, m_0 > d(d-1) .
\ee
Then there exists a global-in-time weak solution $f \in C([0,\infty); \mc{P}(\bt^d \times \br^d))$ of the VPME system \eqref{eq:vpme} with initial data $f_0$. This is the unique solution of \eqref{eq:vpme} with initial datum $f_0$ such that 
$$
\rho_f \in L^\infty_{\text{loc}}([0,+\infty) ; L^\infty(\TT^d) ).
$$
In addition, if $f_0$ has compact support, then at each time $t$, $f(t)$ has compact support.
\end{thm}

This theorem asks for no regularity on $f_0$, only that $f_0 \in L^1 \cap L^\infty(\bt^d \times \br^d)$. The resulting solutions are therefore not $C^1$ classical solutions in general. It is thus useful to introduce a concept of {\it strong} solutions:  the class of bounded distributional solutions $f$ of \eqref{eq:vpme} whose density $\rho_f$ is uniformly bounded: $\rho_f \in L^\infty_{\text{loc}}([0,+\infty) ; L^\infty(\mc{X}) )$. Strong solutions have several convenient properties: in particular, their characteristic ODE system is well-posed and the resulting flow can be used to represent the solutions. A consequence of this is that if the initial datum $f_0$ is additionally assumed to be $C^1$, then the resulting strong solution is in fact a $C^1$ classical solution. Therefore we may also deduce global well-posedness for classical solutions of the VPME system.

In a forthcoming paper, we also consider the problem posed on the whole space; we are able to prove the following global well-posedness result for the whole space systems \eqref{eq:vpme-whole-variable}and \eqref{eq:vpme-whole-fixed}.

\begin{thm}[Global well-posedness: $\RR^3$]
\label{thm:wp}
Let ${f_0 \in L^1\cap L^\infty(\br^3 \times \br^3)}$ be a probability density 
satisfy
$$
f_0(x,v)\leq \frac{C}{(1+|v|)^r} \,\,\,\text{ for some $r>3$}, \qquad \int_{\br^3\times \br^3}|v|^{m_0}f_0(x,v)\,dx\,dv <+\infty \,\,\,\text{ for some $m_0>6$} .
$$
Assume that $g\in L^1\cap L^\infty(\br^3)$, with $g\geq 0$ satisfying $\int_{\br^3}g=1$. Then there exists a unique solution ${f\in L^\infty([0,T] ; L^1\cap L^\infty(\br^3 \times \br^3))}$ of  \eqref{eq:vpme-whole-variable}  (resp. \eqref{eq:vpme-whole-fixed}) with initial datum $f_0$ such that $\rho_f \in L^{\infty}([0,T] ; L^\infty(\br^3))$.
\end{thm}

\begin{remark}
In particular, these results provide well-posedness for the VPME system under the same conditions as were previously known for the Vlasov-Poisson system.
\end{remark}

\subsection{Strategy for $\TT^d$}

\subsubsection{Analysis of the Electric Field}

The first step of the proof is to obtain estimates on the regularity of the electric field $E$. We begin with a decomposition of the electric field, as was used in \cite{IHK1} for the one dimensional setting. The electric field $E$ can be seen as a sum of the electric field appearing in the electron model \eqref{eq:VP}, plus a more regular nonlinear term. For this, we use the notation $E = \bar E +\widehat E $, where
\be
\bar E =-\nabla \bar U ,\qquad \widehat E =-\nabla \widehat U ,
\ee
and $\bar U $ and $\widehat U $ solve respectively
\be \label{electric-field-strategy}
 \Delta \bar U =1-\rho_f  ,\qquad  \Delta \widehat U =e^{\bar U +\widehat U } - 1.
\ee
We expect $\widehat E$ to be more regular than $\bar E$. The key point is to prove this rigorously, taking into account the nonlinearity in the equation satisfied by $\widehat U$. In particular we need to quantify the gain of regularity carefully.

To analyse $\widehat E$, we use techniques from the calculus of variations which allow us to deal with the nonlinearity in the equation for $\widehat U$. We then wish to quantify the gain of regularity in terms of its dependence on $\rho_f$. The key lemma is the following regularity estimate.

\begin{lem} \label{lem:hatU-reg}
Let $d=2,3$.
Assume that $\rho_f \in L^{\frac{d+2}{d}}$. There exist unique $\bar U, \widehat U \in W^{1,2}(\TT^d)$ such that
\be
 \Delta \bar U =1-\rho_f  ,\qquad  \Delta \widehat U =e^{\bar U +\widehat U } - 1.
\ee
Moreover, there exists $\alpha > 0$ such that $\widehat U \in C^{2,\alpha}(\TT^d)$, with the quantitative estimate
\be
\|\widehat U  \|_{C^{2,\alpha}(\bt^d)}  \leq C_{\alpha,d}\,\exp\,\exp  {\Bigl(C_{\alpha,d}\,    \Bigl(1 +  \lVert \rho_f \rVert_{L^{\frac{d+2}{d}}(\bt^d)} \Bigr) \Bigr)}, \qquad\alpha \in \begin{cases} (0,1) \text{ if } d=2 \\ (0, \frac{1}{5}]  \text{ if } d=3. \end{cases} 
\ee
\end{lem}

The choice of $(d+2)/d$ as the integrability exponent is relevant because
this is a quantity that we expect to be bounded uniformly in time, as a consequence of the conservation of the following energy functional associated to the VPME system:
\be \label{def:Ee}
\mc{E} [f ] := \frac{1}{2}\int_{\bt^d \times \br^d} |v|^2 f  \di x \di v + \frac{1 }{2} \int_{\bt^d} |\nabla U |^2 \di x +  \int_{\bt^d} U  e^{U } \di x .
\ee

\begin{lem} \label{lem:rho-Lp}
Let $f \geq 0$ satisfy, for some constant $C_0 > 0$,
\be
\lVert f \rVert_{L^{\infty}(\bt^d \times \br^d)} \leq C_0,\qquad
\mc{E} [f] \leq C_0 ,
\ee
where $\mc{E}$ is the energy functional defined in \eqref{def:Ee}. Then the mass density
\be \label{def:rho}
\rho_f(x) : = \int_{\br^d} f(x,v) \di v
\ee
lies in $L^{(d+2)/d}(\bt^d)$ with
\be \label{rho-Lp}
\lVert \rho_f \rVert_{L^{\frac{d+2}{d}}(\bt^d)} \leq C_1 .
\ee
for some constant $C_1 > 0$ depending on $C_0$ and $d$ only.
\end{lem}

Using these estimates on the electric field, the proof of well-posedness is carried out in two main steps. First we prove the uniqueness of solutions for VPME under the condition that the mass density $\rho_f$ is bounded in $L^\infty(\TT^d)$. Then, we show the global existence of solutions with bounded density, given the assumptions of Theorem~\ref{thm:main}.

\subsubsection{Uniqueness}

The first part of the proof of well-posedness is to prove the uniqueness of strong solutions, i.e. uniqueness under the condition that
\be
\rho_f \in L^\infty_{\text{loc}} \left ( [0, + \infty) ; L^\infty(\TT^d) \right ).
\ee 
For the electron Vlasov-Poisson system \eqref{eq:VP}, Loeper \cite{Loeper} proved uniqueness of solutions under this condition. In the VPME setting, we make use of Loeper's strategy to handle the electric field $\bar E$. However, to deal with $\widehat E$ further nontrivial estimates are necessary. We prove the following estimate, which quantifies the stability of $\widehat E$ with respect to the charge density $\rho_f$. 

\begin{lem} \label{prop:Ustab}
For each $i=1,2$, let $\bar U _i$ and $\widehat U_i$ be respectively solutions of
\be 
 \Delta \bar U _i= h_i - 1, \qquad  \Delta \widehat U _i= e^{\bar U _i + \widehat U _i} - 1 .
\ee
where $h_i \in L^{\infty} \cap L^{(d+2)/d}(\bt^d)$. Then there exists a constant $C_d > 0$ such that
\begin{align}
\label{stab-Uhat}
\lVert \nabla \widehat U _1 - \nabla \widehat U _2 \rVert^2_{L^2(\bt^d)} &\leq \exp\,\exp  {\left [ C_d    \left ( 1 + \max_i\, \lVert h_i \rVert_{L^{(d+2)/d}(\bt^d)}  \right ) \right ]}
 \max_i\, \lVert h_i \rVert_{L^{\infty}(\bt^d)} \, W^2_2(h_1, h_2).
\end{align}
\end{lem}

Using these estimates, we are able to prove the following stability estimate for solutions of the VPME system \eqref{eq:vpme} relative to the initial datum, quantified in the second order Wasserstein distance $W_2$. Uniqueness of strong solutions then follows immediately.

\begin{prop}[Stability for solutions with bounded density] \label{prop:Wstab}
For $i = 1,2$, let $f _i$ be solutions of \eqref{eq:vpme} satisfying for some constant $M$ and all $t \in [0,T]$,
\be \label{str-str_rho-hyp}
\rho[f _i(t)] \leq M .
\ee
Then there exists a constant $C $, depending on $M$, such that, for all $t \in [0,T]$,
\be
W_2\(f _1(t), f _2(t)\)^2 \leq \begin{cases}
16 d e \exp{\biggl [ \log{\frac{W_2\(f _1(0), \, f _2(0)\)^2}{16 d e}}  e^{-C  t} \biggr ]} & \text{ if } t \leq t_0 \\
\max \Big \{W_2\(f _1(0), f _2(0)\)^2 , d \Big \} \, e^{C (1 + \log{16}) (t - t_0)} & \text{ if } t > t_0 .
\end{cases}
\ee
where the time $t_0$ is defined by
\be
t_0 = t_0\big(W_2\(f _1(0), f _2(0)\) \big) = \inf \left\{ t \geq 0 : 16 d e \exp{\biggl [ \log{\frac{W_2\(f _1(0), \, f _2(0)\)^2}{16 d e}}  e^{-C  t} \biggr ]}> d\right\} .
\ee

\end{prop}

\subsubsection{Existence of Solutions}

The proof of existence is based on controlling the moments of solutions. We first show an a priori estimate, proving that the VPME propagates velocity moments of sufficiently high order. This approach was previously used to prove global existence for the electron Vlasov-Poisson system, going back to the work of Lions and Perthame \cite{Lions-Perthame} for the problem posed on $\RR^3$. Pallard \cite{Pallard} proved propagation of moments on the torus and extended the range of moments that could be propagated in the whole space, while Chen and Chen \cite{Chen-Chen} further extended the range of moments available for the torus case.
By extending these methods to the VPME case, we show global-in-time existence of solutions for the VPME system, for any initial datum $f_0 \in L^1 \cap L^\infty(\TT^d \times \RR^d)$ that has a finite velocity moment of order $m_0 > d$. Note that Theorem~\ref{thm:main} requires moments of higher order than this, for the reason that stronger assumptions are required to show uniqueness.

The proposition below shows the propagation of moments for classical solutions of the VPME system. In the proof, the estimates from Lemma~\ref{lem:hatU-reg} on the nonlinear part of the potential $\widehat U$ are crucial.

\begin{prop} \label{prop:moment-propagation}
Let the dimension $d=2$ or $d=3$. Let $0 \leq f_0 \in L^1 \cap L^\infty(\TT^d \times \RR^d)$ have a finite energy and finite velocity moment of order $m_0 > d$:
\be
\mc{E} [f] \leq C_0<+\infty,\qquad \int_{\TT^d \times \RR^d} |v|^{m_0} f_0(x,v) \di x \di v =M_0< +\infty.
\ee
Let $f$ be a $C^1$ compactly supported solution of the VPME system \eqref{eq:vpme}.
Then, for all $T > 0$,
\be
\sup_{[0,T]} \int_{\TT^d \times \RR^d} |v|^{m_0} f(t,x,v) \di x \di v \le C(T, C_0,M_0,m_0,\|f_0\|_\infty).
\ee
\end{prop}

Using this estimate, we then prove the global existence of solutions for the VPME system under these assumptions.We first consider a regularized version of the VPME system:
\be \label{eq:vpme-reg}
\begin{cases}
\partial_t f + v \cdot \nabla_x f - \chi_r \ast_x \nabla_x U \cdot \nabla_v f = 0, \\
 \Delta U = e^U - \chi_r \ast_x \rho_f , \\ \ds
f \vert_{t=0} = f_0 , \; \int_{\TT^d \times \RR^d} f_0 \di x \di v = 1.
\end{cases}
\ee
Here $\chi_r$ is a mollifier defined for $r > 0$ by
\be
\chi_r(x) : = r^{-d} \chi \left ( \frac{x}{r} \right ), \qquad \chi \in C^\infty_c(\bt^d ; [0, +\infty)),
\ee
where $\chi$ is a fixed smooth, radially symmetric function with compact support. 

The regularized system \eqref{eq:vpme-reg} is globally well-posed. This can be proved using standard methods, for example by adapting the approach of Dobrushin \cite{Dobrushin}. The proof of Proposition~\ref{prop:moment-propagation} then provides moment estimates for the solutions of \eqref{eq:vpme-reg} that are uniform in the regularization parameter. We can then extract a limit point and show that it is a global solution of the VPME system. With this method of construction, no regularity is required on the initial datum $f_0$. Moreover, the conservation of the energy $\mc{E}[f]$ defined in \eqref{def:Ee} also follows -- in comparison, the energy of the weak solutions constructed by Bouchut \cite{Bouchut} is non-increasing but not necessarily conserved. We obtain the following existence result.

\begin{thm} \label{thm:existence}
Let $d = 2, 3$. Consider an initial datum $f_0 \in L^1 \cap L^\infty(\bt^d \times \br^d)$ satisfying
\be
\int_{\TT^d \times \RR^d} |v|^{m_0} f_0(x,v) \di x \di v < + \infty, \; \;\mbox{for some}\,\, m_0 > d .
\ee
Then there exists a global-in-time weak solution $f \in C([0,\infty); \mc{P}(\bt^d \times \br^d))$ of the VPME system \eqref{eq:vpme} with initial data $f_0$, such that for all $T>0$,
\be
\sup_{t \in [0,T]} \int_{\TT^d \times \RR^d} |v|^{m_0} f(t,x,v) \di x \di v < + \infty .
\ee
\end{thm}

The proof of Theorem~\ref{thm:main} is then completed by showing that, under the specified decay and moment assumption on $f_0$, the solution provided by Theorem~\ref{thm:existence} has bounded density.
Proposition~\ref{prop:Wstab} then applies, proving the uniqueness of this solution.

\subsection{Strategy for $\RR^3$}

In the whole space case, the overall strategy is similar to the torus case: we first analyse the electrostatic potential using the decomposition $U = \bar U + \widehat U$, where 
\be \label{eq:barU}
- \Delta \bar U = \rho_f, \qquad \lim_{|x| \to 0} \bar U(x) = 0 ,
\ee
and the remainder $\widehat U$ satisfies either
\be \label{eq:hatU-both}
\Delta \widehat U = g e^{\bar U + \widehat U} \quad \text{or} \quad \Delta \widehat U = \frac{ g e^{\bar U + \widehat U}}{\int_{\br^3} g e^{\bar U + \widehat U} \di x}.
\ee
Once again, by using techniques from the calculus of variations we can show that the nonlinear remainder $\widehat U$ is more regular than $\bar U$. However, one first difference with the torus case is that we have to account for the behaviour of the potential at infinity.

A more significant difference occurs for the fixed charge model. Due to the normalisation of the electron charge, the nonlinearity takes a different form compared to the torus case. To deal with this, we use a different functional in the calculus of variations approach to the analysis of $\widehat U$.

For the uniqueness of strong solutions, once again we prove a stability estimate in $W_2$ using stability estimates for the electric field with respect to the charge density $\rho_f$. For $\bar E$ we use estimates devised by Loeper \cite{Loeper}. For $\widehat E$ we again need a version of Lemma~\ref{prop:Ustab}, modified in the fixed charge case to handle the different nonlinearity.

To prove existence, we again use the propagation of moments. However the proof of the propagation of moments in the whole space is very different with respect to the propagation of moments on the torus, and we rely on the approach of Lions and Perthame \cite{Lions-Perthame}, making use of the regularity estimates on $\widehat U$.

\bibliography{proc-VPME-bib}

\begin{thebibliography}{10}

\bibitem{Arsenev}
A.~Arsenev.
\newblock Existence in the large of a weak solution to the {Vlasov} system of
  equations.
\newblock {\em Zh. Vychisl. Mat. i Mat. Fiz.}, 15:136--147, 1975.

\bibitem{BGNS18}
C.~Bardos, F.~Golse, T.~T. Nguyen, and R.~Sentis.
\newblock The {M}axwell-{B}oltzmann approximation for ion kinetic modeling.
\newblock {\em Phys. D}, 376/377:94--107, 2018.

\bibitem{Batt-Rein}
J.~Batt and G.~Rein.
\newblock {Global classical solutions of the periodic Vlasov-Poisson system in
  three dimensions}.
\newblock {\em C. R. Acad. Sci. Paris S{\'{e}}r. I Math.}, 313(6):411--416,
  1991.

\bibitem{Bellan}
P.~M. Bellan.
\newblock {\em Fundamentals of Plasma Physics}.
\newblock Cambridge University Press, 2008.

\bibitem{BPLT1991}
G.~Bonhomme, T.~Pierre, G.~Leclert, and J.~Trulsen.
\newblock Ion phase space vortices in ion beam-plasma systems and their
  relation with the ion acoustic instability: numerical and experimental
  results.
\newblock {\em Plasma Physics and Controlled Fusion}, 33(5):507--520, may 1991.

\bibitem{Bouchut}
F.~Bouchut.
\newblock {Global weak solution of the Vlasov-Poisson system for small
  electrons mass}.
\newblock {\em Comm. Partial Differential Equations}, 16(8-9):1337--1365, 1991.

\bibitem{Bouchut-Dolbeault95}
F.~Bouchut and J.~Dolbeault.
\newblock On long time asymptotics of the {Vlasov-Fokker-Planck} equation and
  of the {Vlasov-Poisson-Fokker-Planck} system with {Coulombic} and {Newtonian}
  potentials.
\newblock {\em Differential Integral Equations}, 8(3):487--514, 1995.

\bibitem{Braun-Hepp}
W.~Braun and K.~Hepp.
\newblock {The Vlasov dynamics and its fluctuations in the {$1/N$} limit of
  interacting classical particles}.
\newblock {\em Comm. Math. Phys.}, 56(2):101--113, 1977.

\bibitem{Chen-Chen}
Z.~Chen and J.~Chen.
\newblock Moments propagation for weak solutions of the {Vlasov-Poisson} system
  in the three-dimensional torus.
\newblock {\em J. Math. Anal. Appl.}, 42(1):728--737, 2019.

\bibitem{Desvillettes-Villani2000}
L.~Desvillettes and C.~Villani.
\newblock On the spatially homogeneous {Landau} equation for hard potentials
  part ii : h-theorem and applications.
\newblock {\em Communications in Partial Differential Equations},
  25(1-2):261--298, 2000.

\bibitem{Desvillettes-Villani2005}
L.~Desvillettes and C.~Villani.
\newblock On the trend to global equilibrium for spatially inhomogeneous
  kinetic systems: {The} {Boltzmann} equation.
\newblock {\em Invent. Math.}, 159:245--316, 2005.

\bibitem{Dobrushin}
R.~L. Dobrushin.
\newblock {Vlasov Equations}.
\newblock {\em Funktsional. Anal. i Prilozhen.}, 13(2):48--58, 1979.

\bibitem{IGP-WP}
M.~Griffin-Pickering and M.~Iacobelli.
\newblock Global well-posedness in 3-dimensions for the {Vlasov--Poisson}
  system with massless electrons.
\newblock arXiv:1810.06928.

\bibitem{Gurevich-Pitaevsky75}
A.~V. Gurevich and L.~P. Pitaevsky.
\newblock Non-linear dynamics of a rarefied ionized gas.
\newblock {\em Progress in Aerospace Sciences}, 16(3):227 -- 272, 1975.

\bibitem{IHK1}
D.~Han-Kwan and M.~Iacobelli.
\newblock {The quasineutral limit of the Vlasov-Poisson equation in Wasserstein
  metric}.
\newblock {\em Commun. Math. Sci.}, 15(2):481--509, 2 2017.

\bibitem{Hauray14}
M.~Hauray.
\newblock Mean field limit for the one dimensional {V}lasov-{P}oisson equation.
\newblock In {\em S\'{e}minaire {L}aurent {S}chwartz---\'{E}quations aux
  d\'{e}riv\'{e}es partielles et applications. {A}nn\'{e}e 2012--2013}, Exp.
  No. XXI, S\'{e}min. \'{E}qu. D\'{e}riv. Partielles. \'{E}cole Polytech.,
  Palaiseau, 2014.

\bibitem{Herda16}
M.~Herda.
\newblock On massless electron limit for a multispecies kinetic system with
  external magnetic field.
\newblock {\em Journal of Differential Equations}, 260(11):7861 -- 7891, 2016.

\bibitem{Horst1993}
E.~Horst.
\newblock On the asymptotic growth of the solutions of the {V}lasov-{P}oisson
  system.
\newblock {\em Math. Methods Appl. Sci.}, 16(2):75--86, 1993.

\bibitem{Horst-Hunze}
E.~Horst and R.~Hunze.
\newblock Weak solutions of the initial value problem for the unmodified non-
  linear {Vlasov} equation.
\newblock {\em Math. Methods Appl. Sci.}, 6(2):262--279, 1984.

\bibitem{Jeans}
J.~H. Jeans.
\newblock On the theory of star-streaming and the structure of the universe.
\newblock {\em Monthly Notices of the Royal Astronomical Society}, 76:70--84,
  1915.

\bibitem{Lifshitz-Pitaevskii}
E.~M. Lifshitz and L.~P. Pitaevskii.
\newblock {\em Physical Kinetics}, volume~10 of {\em Course of Theoretical
  Physics}.
\newblock Pergamon Press, 1981.

\bibitem{Lions-Perthame}
P.~L. Lions and B.~Perthame.
\newblock {Propagation of moments and regularity for the 3-dimensional
  Vlasov-Poisson system}.
\newblock {\em Invent. Math.}, 105(2):415--430, 1991.

\bibitem{Loeper}
G.~Loeper.
\newblock {Uniqueness of the solution to the Vlasov-Poisson system with bounded
  density}.
\newblock {\em J. Math. Pures Appl. (9)}, 86(1):68--79, 2006.

\bibitem{Mason71}
R.~J. Mason.
\newblock Computer simulation of ion-acoustic shocks. {T}he diaphragm problem.
\newblock {\em The Physics of Fluids}, 14(9):1943--1958, 1971.

\bibitem{Medvedev2011}
Y.~V. Medvedev.
\newblock Ion front in an expanding collisionless plasma.
\newblock {\em Plasma Physics and Controlled Fusion}, 53(12):125007, nov 2011.

\bibitem{Pallard}
C.~Pallard.
\newblock Moment propagation for weak solutions to the {Vlasov-Poisson} system.
\newblock {\em Comm. Partial Differential Equations}, 37(7):1273--1285, 2012.

\bibitem{Pfaffelmoser}
K.~Pfaffelmoser.
\newblock {Global classical solutions of the Vlasov-Poisson system in three
  dimensions for general initial data}.
\newblock {\em J. Differential Equations}, 95(2):281--303, 1992.

\bibitem{Robert}
R.~Robert.
\newblock Unicit{\'{e}} de la solution faible {\`{a}} support compact de
  l'{\'{e}}quation de {Vlasov--Poisson}.
\newblock {\em C. R. Acad. Sci. Paris S{\'{e}}r. I Math.}, 324(8):873--877,
  1997.

\bibitem{SCM}
P.~Sakanaka, C.~Chu, and T.~Marshall.
\newblock Formation of ion-acoustic collisionless shocks.
\newblock {\em The Physics of Fluids}, 14(611), 1971.

\bibitem{Schaeffer}
J.~Schaeffer.
\newblock {Global existence of smooth solutions to the Vlasov-Poisson system in
  three dimensions}.
\newblock {\em Comm. Partial Differential Equations}, 16(8-9):1313--1335, 1991.

\bibitem{Ukai-Okabe}
S.~Ukai and T.~Okabe.
\newblock On classical solutions in the large in time of two-dimensional
  {Vlasov's} equation.
\newblock {\em Osaka J. Math.}, 15(2):245--261, 1978.

\bibitem{Villani1996}
C.~Villani.
\newblock On the {Cauchy} problem for {Landau} equation: sequential stability,
  global existence.
\newblock {\em Adv. Differential Equations}, 1(5):793--816, 1996.

\bibitem{Vlasov1}
A.~A. Vlasov.
\newblock On the vibration properties of the electron gas.
\newblock {\em Zh. Eksper. Teor. Fiz.}, 8(3):291, 1938.

\end{thebibliography}
\bibliographystyle{abbrv}

\end{document}